# SURFACES FAMILY WITH COMMON NULL ASYMPTOTIC

[1] Gülnur Şaffak Atalay, [2] Emin Kasap

[1] e-mail: gulnur.saffak@omu.edu.tr

[2] e-mail: kasape@omu.edu.tr

Department of Mathematics, Arts and Science Faculty,

Ondokuz Mayis University, Samsun 55139, Turkey

**Abstract**

We analyzed the problem of finding a surfaces family through an asymptotic curve with Cartan frame. We obtain the parametric representation for surfaces family whose members have the same as an asymptotic curve. By using the Cartan frame of the given null curve, we present the surface as a linear combination of this frame and analysed the necessary and sufficient condition for that curve to satisfy the asymptotic requirement. We illustrate the method by giving some examples.

## 1.Introduction

Asymptotic curve on a surface has been a long-term research topic in Differential Geometry [1,2,3]. A useful interpretation of asymptotic directions and asymptotic curves is given by the Dupin indicatrix [4]. Lie et al.[5] derived the necessary an sufficient condition for a given curve to be the line of curvature on the surface. Wang et al.[6] studied the problem of costructing a surface family from a given spatial geodesic. Kasap et. al. [7] generalized the work of Wang by introducing new types of marching-scale functions, coefficients of the Frenet frame appearing in the parametric representation of surfaces. In [8] Kasap and Akyıldız defined surfaces with a common geodesic in Minkowski 3-space and gave the sufficient conditions on marching-scale functions so that the given curve is a common on that surfaces. Şaffak and Kasap [9] studied family of surfaces with a common null geodesic. Recently, in [10] Bayram et. al. defined the surface pencil with a common asymptotic curve. They introduced three types of marching-scale functions and derived the necessary and sufficient conditions on them to satisfy both parametric and asymptotic requirements.

In this paper, we first take a null curve $\alpha = \alpha(s)$ on surface $\varphi = \varphi(s,t)$. Using Cartan frame of the curve, we derive the necessary and sufficient conditions such that the curve is an asymptotic on the surface $\varphi = \varphi(s,t)$. Thus, we define the family of surfaces with common null asymptotic. Finally, we demonstrated some interesting surfaces about subject.

## 2. Preliminaries

Let $IR_1^3$ be Minkowski 3-space with natural Lorentzian metric $<.,.> = -dx^2 + dy^2 + dz^2$
A Cartan frame denoted by $\{\ell(s), n(s), u(s)\}$ on a null curve $\alpha = \alpha(s)$ with tangent vector $\ell = \ell(s)$ where $<n,n> = 0$, $<u,u> = 1$, $<\ell,n> = -1$, $<\ell,u> = 0$, $<n,u> = 1$ and $\det(\ell, n, u) < 0$ ,[11].

The Cartan frame $\{\ell(s), n(s), u(s)\}$ satisfies the following Frenet -Serret formula:

$$\ell' = k_1 u, \quad n' = -k_2 u, \quad u' = -k_2 \ell + k_1 n \tag{2.1}$$

where $k_1 = \|\ell'\|$ and $k_2 = -<n', u>$ are are called the curvature functions of $\alpha = \alpha(s)$ [11].

For the Cartan frame $\{\ell(s), n(s), u(s)\}$, it is easy to see that

$$\ell \times n = u, \quad n \times u = -n, \quad \ell \times u = \ell \tag{2.2}$$

A parametric curve $\alpha(s)$ is a curve on a surface $\varphi = \varphi(s,t)$ in $IR_1^3$ that has a constant s or t-parameter value. In other words, there exists a parameter $s_o$ or $t_o$ such that $\alpha(t) = \varphi(s_o, t)$ or $\alpha(s) = \varphi(s, t_o)$.

## 3. Surfaces with Common Null Asymptotic

$\alpha = \alpha(s)$ be an parametric null curve with Cartan frame $\{\ell, n, u\}$ on surface $\varphi = \varphi(s,t)$. We assume that $\|\alpha''(s)\| \neq 0$.

The surface $\varphi$ is defined by

$$\varphi(s,t) = \alpha(s) + [x(s,t)\ell(s) + y(s,t)n(s) + z(s,t)u(s)] \tag{3.1}$$

where $x(s,t)$, $y(s,t)$ and $z(s,t)$ are $C^1$ functions. The values of the marching-scale functions $x(s,t)$, $y(s,t)$ and $z(s,t)$ indicate, respectively, the extension-like, flexion-like and retortion-like effects, by the point unit through the time t, starting from $\alpha(s)$.

Our goal is to find the necessary and sufficient conditions for which the null curve $\alpha(s)$ is a parameter and an asymptotic curve on the surface $\varphi(s,t)$.

Firstly, since $\alpha(s)$ is an parametric curve on surface $\varphi(s,t)$, there exists a parameter $t_0 \in [T_1, T_2]$ such that

$$x(s,t_o) = y(s,t_o) = z(s,t_o) = 0 \qquad L_1 \leq s \leq L_2, \; T_1 \leq t_0 \leq T_2 \tag{3.2}$$

Secondly, according to the asymptotic theory [11], the curve $\alpha(s)$ is an asymptotic curve on the surface $\varphi(s,t)$ if and only if

$$\frac{\partial N}{\partial s}(s,t_0) \cdot l(s) = 0 \tag{3.3}$$

is satisfied, where '.' denotes the Lorentzian inner product.

The normal vector of $\varphi = \varphi(s,t)$ can be written as

$$N(s,t) = \frac{\partial \varphi(s,t)}{\partial s} \times \frac{\partial \varphi(s,t)}{\partial t}.$$

From (2.1) and (3.1), the normal vector can be expressed as

$$N(s,t) = \left[\left(1 + \frac{\partial x(s,t)}{\partial s} - z(s,t)k_2\right)\frac{\partial z(s,t)}{\partial t} - \left(\frac{\partial z(s,t)}{\partial s} + x(s,t)k_1 - y(s,t)k_2\right)\frac{\partial x(s,t)}{\partial t}\right]\ell +$$
$$\left[\left(\frac{\partial z(s,t)}{\partial s} + x(s,t)k_1 - y(s,t)k_2\right)\frac{\partial y(s,t)}{\partial t} - \left(\frac{\partial y(s,t)}{\partial s} + z(s,t)k_1\right)\frac{\partial z(s,t)}{\partial t}\right]n +$$
$$\left[\left(1 + \frac{\partial x(s,t)}{\partial s} - z(s,t)k_2\right)\frac{\partial y(s,t)}{\partial t} - \left(\frac{\partial y(s,t)}{\partial s} + z(s,t)k_1\right)\frac{\partial x(s,t)}{\partial t}\right]u.$$

Thus,

$$N(s,t_0) = \phi_1(s,t_o)l(s) + \phi_2(s,t_o)n(s) + \phi_3(s,t_o)u(s)$$

where

$$\phi_1(s,t_o) = \left(1 + \frac{\partial x(s,t_o)}{\partial s}\right)\frac{\partial z(s,t_o)}{\partial t} - \frac{\partial z(s,t_o)}{\partial s}\frac{\partial x(s,t_o)}{\partial t},$$

$$\phi_2(s,t_o) = \frac{\partial z(s,t_o)}{\partial s}\frac{\partial y(s,t_o)}{\partial t} - \frac{\partial y(s,t_o)}{\partial s}\frac{\partial z(s,t_o)}{\partial t},$$

$$\phi_3(s,t_o) = \left(1 + \frac{\partial x(s,t_o)}{\partial s}\right)\frac{\partial y(s,t_o)}{\partial t} - \frac{\partial y(s,t_o)}{\partial s}\frac{\partial x(s,t_o)}{\partial t}.$$

From the Eqn.(3.3), we should have

$$\frac{\partial N}{\partial s}(s,t_0) \cdot l(s) = 0 \Leftrightarrow \frac{\partial(\phi_1(s,t_o)l(s) + \phi_2(s,t_o)n(s) + \phi_3(s,t_o)u(s))}{\partial s} \cdot l(s) = 0$$

$$\Leftrightarrow \frac{\partial \phi_2}{\partial s}(s,t_o) + k_1 \phi_3(s,t_o) = 0$$

Since $k_1 = \|\alpha''(s)\| \neq 0$, $\frac{\partial \phi_2}{\partial s}(s,t_o) = 0$ and by Eqn. (3.2) we have $\phi_3(s,t_o) = \frac{\partial y(s,t_o)}{\partial t}$.
Therefore, Eqn. (3.3) is simplified to

$$\Leftrightarrow \frac{\partial y(s,t_o)}{\partial t} = 0 \tag{3.4}$$

Firstly, similar with [6], for the purposes of simplification an analysis, we consider the case when the marching-scale functions $x(s,t)$, $y(s,t)$ and $z(s,t)$ are decomposed into two factors;

$$x(s,t) = k(s)X(t),$$
$$y(s,t) = m(s)Y(t), \quad\quad L_1 \leq s \leq L_2, \; T_1 \leq t_0 \leq T_2$$
$$z(s,t) = t(s)Z(t),$$

where $k(s)$, $m(s)$, $t(s)$, $X(t)$, $Y(t)$ and $Z(t)$ are $C'$ functions. Now, we can give the following corollary.

**Corollory 3.1.** The necessary and sufficient condition of the null curve $\alpha(s)$ being an asymptotic line on the surface $\varphi(s,t)$ is

$$\begin{cases} X(t_0) = Y(t_0) = Z(t_0) = 0 \quad (x(s,t_0) = y(s,t_0) = z(s,t_0) = 0) \\ m(s) = 0 \text{ or } \frac{dY(t_0)}{dt} = 0 \end{cases} \tag{3.5}$$

Now let us consider other types of the maching-scale functions. In the Eqn. (3.1), maching-scale functions $x(s,t), y(s,t)$ and $z(s,t)$ can be choosen in two different forms;

1) If we choose;

$$x(s,t) = \sum_{i=1}^{p} a_{1i} k(s)^i x(t)^i$$

$$y(s,t) = \sum_{i=1}^{p} a_{2i} m(s)^i y(t)^i \qquad (3.6)$$

$$z(s,t) = \sum_{i=1}^{p} a_{3i} t(s)^i z(t)^i$$

then we can simply express the sufficient condition for which the curve $\alpha(s)$ is an asymptotic curve on the surface $\varphi(s,t)$ as

$$\begin{cases} X(t_0) = Y(t_0) = Z(t_0) = 0 \quad (x(s,t_0) = y(s,t_0) = z(s,t_0) = 0) \\ a_{21} = 0 \text{ or } m(s) = 0 \text{ or } \dfrac{dY(t_0)}{dt} = 0 \end{cases} \qquad (3.7)$$

where $k(s), m(s), t(s), X(t), Y(t)$ and $Z(t)$ are $C^1$ functions, $a_{ij} \in IR$, $i = 1,2,3$, $j = 1,2,...,p$.

2) If we take

$$x(s,t) = f\left(\sum_{i=1}^{p} a_{1i} k(s)^i x(t)^i\right)$$

$$y(s,t) = g\left(\sum_{i=1}^{p} a_{2i} m(s)^i y(t)^i\right) \qquad (3.8)$$

$$z(s,t) = h\left(\sum_{i=1}^{p} a_{3i} t(s)^i z(t)^i\right)$$

then we can write the sufficient condition for which the curve $\alpha(s)$ is an asymptotic curve on the surface $\varphi(s,t)$ as

$$\begin{cases} X(t_0) = Y(t_0) = Z(t_0) = 0 \quad (x(s,t_0) = y(s,t_0) = z(s,t_0) = 0) \\ a_{21} = 0 \text{ or } m(s) = 0 \text{ or } g'(0) = 0 \text{ or } \dfrac{dY(t_0)}{dt} = 0 \end{cases} \qquad (3.9)$$

where $k(s), m(s), t(s), X(t), Y(t), Z(t), f, g$ and $h$ are $C^1$ functions.

**Example 3.1**

Let $\alpha(s) = (s, \sin(s), \cos(s))$ be a null curve. Then $\alpha$ is framed by

$$\ell(s) = (1, \cos(s), -\sin(s)),$$
$$n(s) = \left(\frac{1}{2}, -\frac{1}{2}\cos(s), \frac{1}{2}\sin(s)\right),$$
$$u(s) = (0, -\sin(s), -\cos(s)).$$

**a)** If we take, $x(s,t) = t$, $y(s,t) = z(s,t) = 0$ and $t_o = 0$ where $0 \leq s \leq 2\pi$ and $-3 \leq t \leq 3$, then Eqn.(3.5) is satisfied. Thus, we obtain member of surfaces family with common null asymptotic $\alpha = \alpha(s)$ as shown in Fig. 1:

$$\varphi_1(s,t) = (s+t, \sin(s)+t\cos(s), \cos(s)-t\sin(s))$$

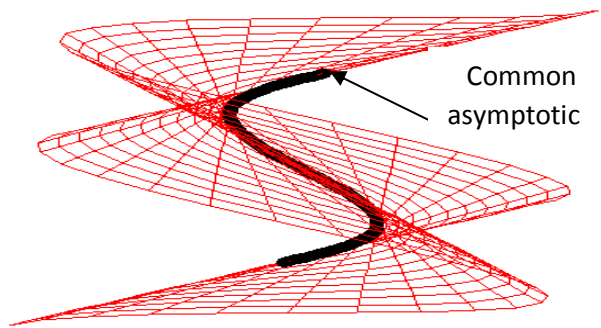

Common asymptotic

Fig.1. A member of surfaces family and its common null asymptotic.

**b)** If we take, $x(s,t) = t^2$, $y(s,t) = 0$, $z(s,t) = t$ and $t_o = 0$ where $0 \leq s \leq 2\pi$ and $0 \leq t \leq 0{\cdot}6$ ,then Eqn.(3.5) is satisfied. Thus, we obtain another member of surfaces family with common null asymptotic $\alpha = \alpha(s)$ as shown in Fig. 2:

$$\varphi_2(s,t) = \left(s+t^2, \sin(s)(1-t)+t^2\cos(s), \cos(s)(1-t)-t^2\sin(s)\right)$$

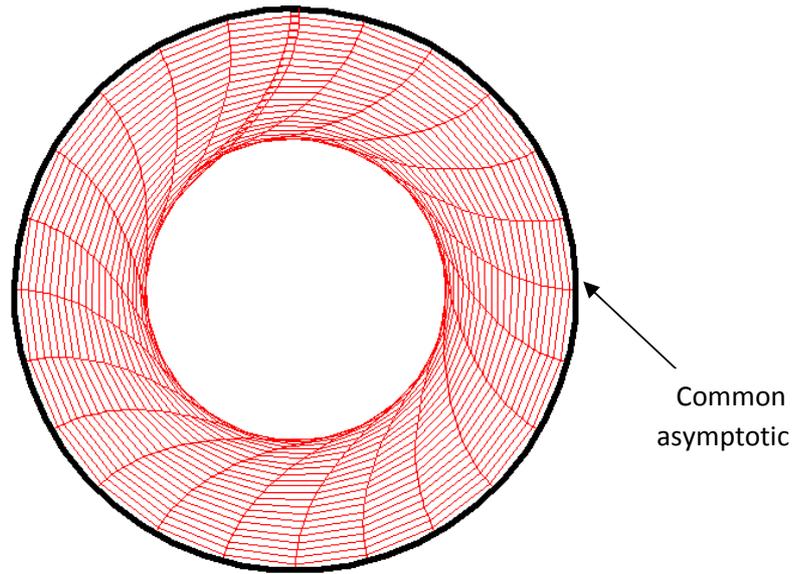

Fig. 2. A member of surfaces family and its common null asymptotic.

**c)** If we take, $x(s,t) = t$, $y(s,t) = 0$, $z(s,t) = t$ and $t_o = 0$ where $0 \leq s \leq 2\pi$ and $0 \leq t \leq 0 \cdot 6$, then Eqn.(3.5) is satisfied. Thus, we obtain another member of surfaces family with common null asymptotic $\alpha = \alpha(s)$ as shown in Fig. 3:

$$\varphi_3(s,t) = \left(s+t, \sin(s)(1-t)+t\cos(s), \cos(s)(1+t)-t\sin(s)\right)$$

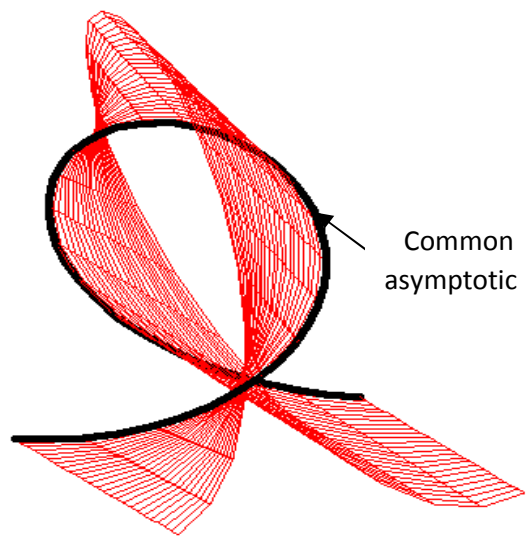

Fig. 3. A member of surfaces family and its common null asymptotic

**d)** If we take, $x(s,t) = st$, $y(s,t) = st^2$, $z(s,t) = t^2$ and $t_o = 0$ where $0 \leq s \leq 6\pi$ and $0 \leq t \leq 0.6$, then Eqn. (3.7) is satisfied. Thus, we obtain another member of surfaces family with common null asymptotic $\alpha = \alpha(s)$ as shown in Fig. 4:

$$\varphi_4(s,t) = \left( s + st + \frac{1}{2}st^2,\ stcos(s) + \sin(s) - \frac{1}{2}st^2 cos(s) - t^2 \sin(s),\ cos(s) - stsin(s) + \frac{1}{2}st^2 sin(s) - t^2 \cos(s) \right)$$

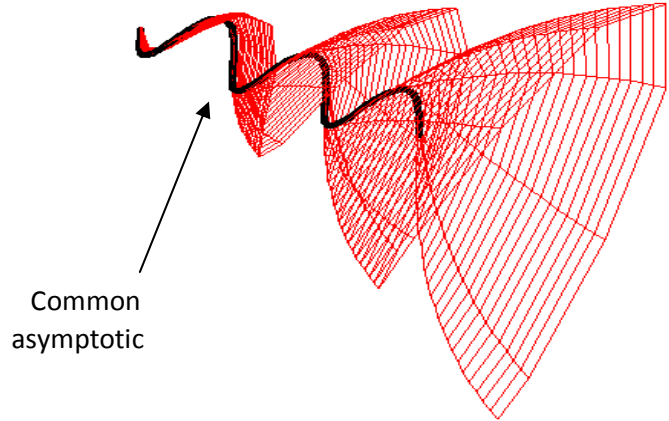

Common asymptotic

Fig. 4. A member of surfaces family and its common null asymptotic.

**Example 3.2.** Let $\alpha(s) = \left( -\frac{\sqrt{2}}{12}s^3 - \frac{\sqrt{2}}{2}s,\ -\frac{s^2}{2},\ -\frac{\sqrt{2}}{12}s^3 + \frac{\sqrt{2}}{2}s \right)$ be a null curve. The curve is called a null cubic. Then $\alpha$ is framed by

$$\ell(s) = \left( -\frac{\sqrt{2}}{4}s^2 - \frac{\sqrt{2}}{2},\ -s,\ -\frac{\sqrt{2}}{4}s^2 + \frac{\sqrt{2}}{2} \right),\quad n(s) = \left( -\frac{\sqrt{2}}{2},\ 0,\ -\frac{\sqrt{2}}{2} \right)\ \text{and}$$

$$u(s) = \left( -\frac{\sqrt{2}}{2}s,\ -1,\ -\frac{\sqrt{2}}{2}s \right).$$

If we take, $x(s,t) = y(s,t) = 0$, $z(s,t) = \frac{2}{\sqrt{2}}t$ and $t_o = 0$ where $-4 \leq s \leq 4$ and $-10 \leq t \leq 10$, then Eqn.(3.5) is satisfied. Thus, we obtain member of surfaces family with common null asymptotic $\alpha = \alpha(s)$ as shown in Fig.5:

$$\varphi_1(s,t) = \left( -\frac{\sqrt{2}}{12}s^3 - \frac{\sqrt{2}}{2}s - ts,\ -\frac{s^2}{2} - \frac{2}{\sqrt{2}}t,\ -\frac{\sqrt{2}}{12}s^3 + \frac{\sqrt{2}}{2}s - ts \right)$$

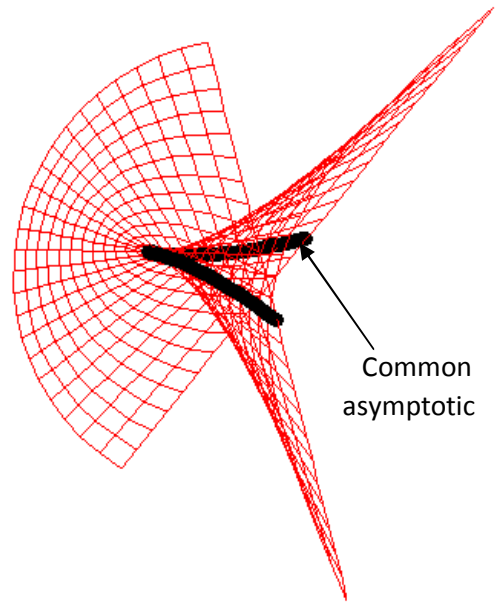

Fig. 5. A member of surfaces family and its common null asymptotic.